\documentclass{amsart}
\usepackage{amsxtra,amssymb,multicol,graphicx}
\newcommand{\oneg}{\vcenter{\hbox{\includegraphics[scale=0.1]{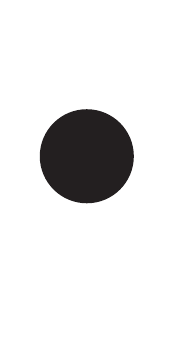}}}}
\newcommand{\twog}{\vcenter{\hbox{\includegraphics[scale=0.1]{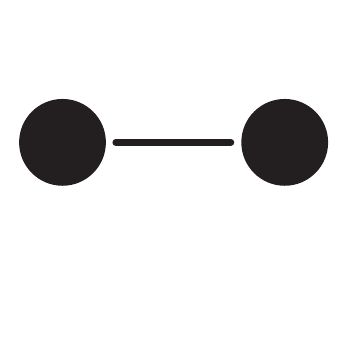}}}}
\newcommand{\threeg}{\vcenter{\hbox{\includegraphics[scale=0.1]{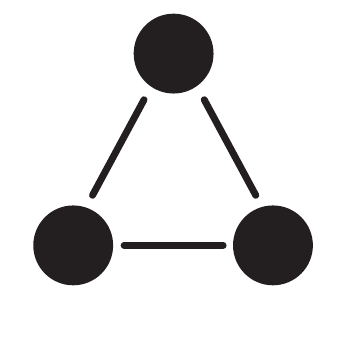}}}}
\newcommand{\fourg}{\vcenter{\hbox{\includegraphics[scale=0.1]{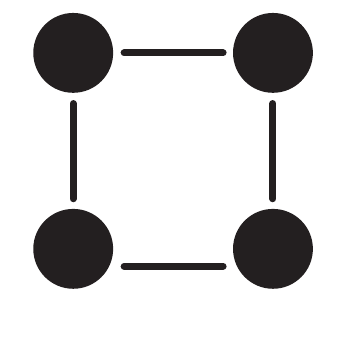}}}}

\title{WHAT IS ...?}
\author{Daniel Glasscock}
\begin{document} 
%\maketitle
\rightline{\Huge \bf WHAT IS ... a Graphon?}
\vskip1cm
\rightline{\it \huge Daniel Glasscock}
\vskip1in

\begin{multicols}{2}

\noindent Large graphs are ubiquitous in mathematics, and describing their structure is an important goal of modern combinatorics. One way to study large, finite objects is to pass from sequences of larger and larger such objects to ideal limiting objects. Done properly, properties of the limiting objects reflect properties of the finite objects which approximate them, and vice versa.

Graphons, short for graph functions, are the limiting objects for sequences of large, finite graphs with respect to the so-called cut metric. They were introduced and developed by C. Borgs, J. T. Chayes, L. Lov\'{a}sz, V. T. S\'{o}s, B. Szegedy, and K. Vesztergombi in \cite{lovasz2008} and \cite{lovasz2006}. Graphons arise naturally wherever sequences of large graphs appear: extremal graph theory, property testing of large graphs, quasi-random graphs, random networks, et cetera.

Let's begin with some definitions and a motivating example. A \emph{graph} $G$ is a set of vertices $V(G)$ and a set of edges $E(G)$ between the vertices (excluding loops and multiple edges). A \emph{graph homomorphism} from $H$ to $G$ is a map from $V(H)$ to $V(G)$ that preserves edge adjacency; that is, for every edge $\{v,w\}$ in $E(H)$, the edge $\{\varphi(v),\varphi(w)\}$ is in $E(G)$.  Denote by $\hom (H,G)$ the number of homomorphisms from $H$ to $G$. For example, $\hom(\oneg,G) = |V(G)|$, $\hom(\twog,G) = 2 |E(G)|$, and $\hom(\threeg,G)$ is 6 times the number of triangles in $G$. Normalizing by the total number of possible maps, we get the \emph{homomorphism density} of $H$ into $G$,
\[t(H,G) = \frac{\hom(H,G)}{|V(G)|^{|V(H)|}},\]
the probability that a randomly chosen map from $V(H)$ to $V(G)$ preserves edge adjacency.  This number also represents the density of $H$ as a subgraph in $G$ asymptotically as $n = |V(G)| \rightarrow \infty$. For example, $t(\twog,G) = 2 |E(G)| / n^2$ while the density of edges in $G$ is $2 |E(G)| / n(n - 1)$; these two expressions are nearly the same when $n$ is large.

Consider the following problem from extremal graph theory:
\begin{center} \emph{How many 4-cycles must there be in a graph with edge density at least $1/2$?} \end{center}
It is easy to see that there are at most on the order of $n^4$ 4-cycles in any graph; a theorem of Erd\H{o}s gives that graphs with at least half the number of possible edges have \emph{at least} on the order of $n^4$ 4-cycles. More specifically, for any graph $G$,
\[t(\fourg,G) \geq t(\twog,G)^4,\]
meaning that if $t(\twog,G) \geq 1/2$, then $t(\fourg,G) \geq 1/16$. In light of this, the problem may be reformulated into a minimization one: \emph{Minimize $t(\fourg,G)$ over finite graphs $G$ satisfying $t(\twog,G) \geq 1/2$}. With some work, it may be shown that no finite graph $G$ with $t(\twog,G) \geq 1/2$ achieves the minimum $t(\fourg,G) = 1/16$.

It's useful at this point to draw an analogy with a problem from elementary analysis: \emph{Minimize $x^3 - 6x$ over rational numbers $x$ satisfying $x \geq 0$}. This polynomial has a unique minimum on $x \geq 0$ at $x = \sqrt 2$, so the best we may do over the rationals is show that the polynomial achieves values approaching this minimum along a sequence of rationals approaching $\sqrt 2$. We know well to avoid this complication by completing the rational numbers to the reals and realizing the limit of such a sequence as $\sqrt 2$.

There is a sequence of finite graphs with edge density at least 1/2 and 4-cycle density approaching 1/16. Let $R_n$ be an instance of a random graph on $n$ vertices where each edge is decided independently with probability $1/2$.  Throwing away those $R_n$'s for which $t(\twog,R_n) < 1/2$, the 4-cycle density in the remaining graph sequence limits to 1/16 almost surely. Following the $\sqrt 2$ analogy, we should look to realize the limit of this sequence of finite graphs and understand how it solves the minimization problem at hand.

What might the limit of the sequence of random graphs $(R_n)_n$ be? From the adjacency matrix of a labeled graph, construct the graph's \emph{pixel picture} by turning the 1's into black squares, erasing the 0's, and scaling to the unit square $[0,1]^2$.

{\centering
\begin{tabular}{ccc}
& & \\
$\left(\rule{0cm}{.8cm}\right.$ \kern-.5em {\small \begin{tabular}{p{.09cm}p{.09cm}p{.09cm}p{.09cm}}
0 & 1 & 0 & 1\\
1 & 0 & 1 & 0 \\
0 & 1 & 0 & 1 \\
1 & 0 & 1 & 0
\end{tabular}} \kern-.5em $\left.\rule{0cm}{.8cm}\right)$ & $\boldsymbol{\longrightarrow}$ & $\vcenter{\hbox{\includegraphics[scale=.39]{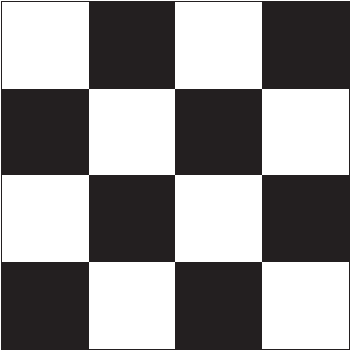}}}$ \\
&  &
\end{tabular}\par}

\noindent Pixel pictures may be seen to ``converge'' graphically; those of larger and larger random graphs with edge probability 1/2, regardless of how they are labeled, seem to converge to a gray square, the constant 1/2 function on $[0,1]^2$. 

\vskip8pt
\begin{center}\includegraphics[scale=.39]{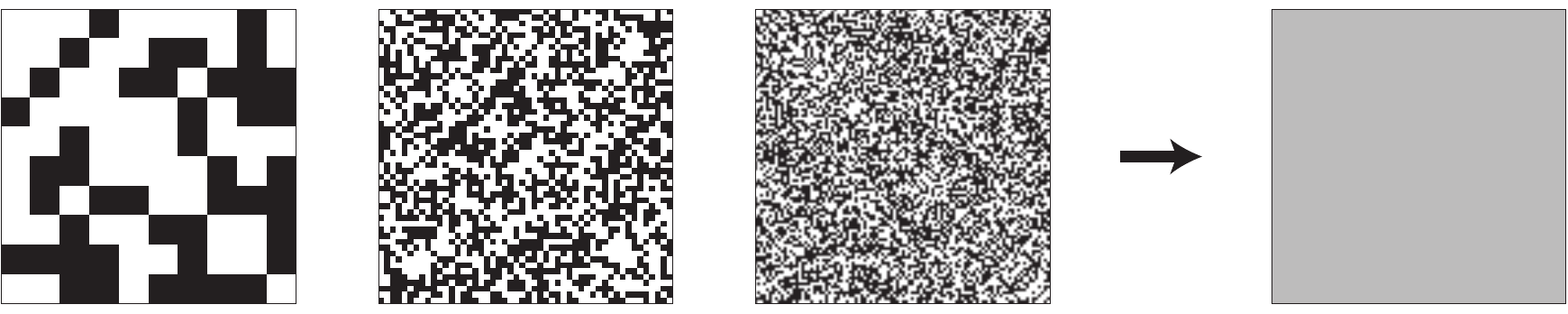}\end{center}
\vskip5pt

The constant 1/2 function on $[0,1]^2$ is an example of a labeled graphon. A \emph{labeled graphon} is a symmetric, Lebesgue-measurable function from $[0,1]^2$ to $[0,1]$ (modulo the usual identification almost everywhere); they may be thought of as edge-weighted graphs on the vertex set $[0,1]$. An \emph{unlabeled graphon} is a graphon up to re-labeling, where a \emph{re-labeling} is the result of applying an invertible, measure preserving transformation to the $[0,1]$ interval. Note that any pixel picture is a labeled graphon, meaning that (labeled) graphs are (labeled) graphons.

As another example of this convergence, consider the \emph{growing uniform attachment} graph sequence $(G_n)_n$ defined inductively as follows. Let $G_1 = \oneg$. For $n \geq 2$, construct $G_n$ from $G_{n-1}$ by adding one new vertex, then, considering each pair of non-adjacent vertices in turn, drawing an edge between them with probability $1/n$.  This sequence almost surely limits to the graphon $1-\max(x,y)$. (Since matrices are indexed with $(0,0)$ in the top left corner, so too are graphons.)

\vskip8pt
\begin{center}\includegraphics[scale=.39]{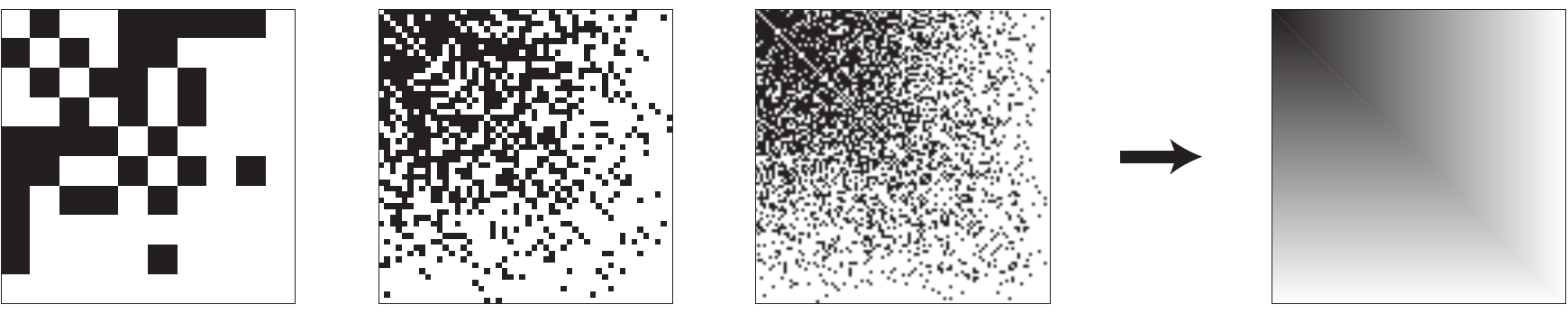}\end{center}
\vskip5pt

There are two natural ways to label a complete bipartite graph, and each suggests a different limit graphon for the complete bipartite graph sequence. Both sequences of labeled graphons in fact have the same limit, as indicated in the diagram; the reader is encouraged to return to this example after we define this convergence more precisely.

\vskip8pt
\begin{center}\includegraphics[scale=.39]{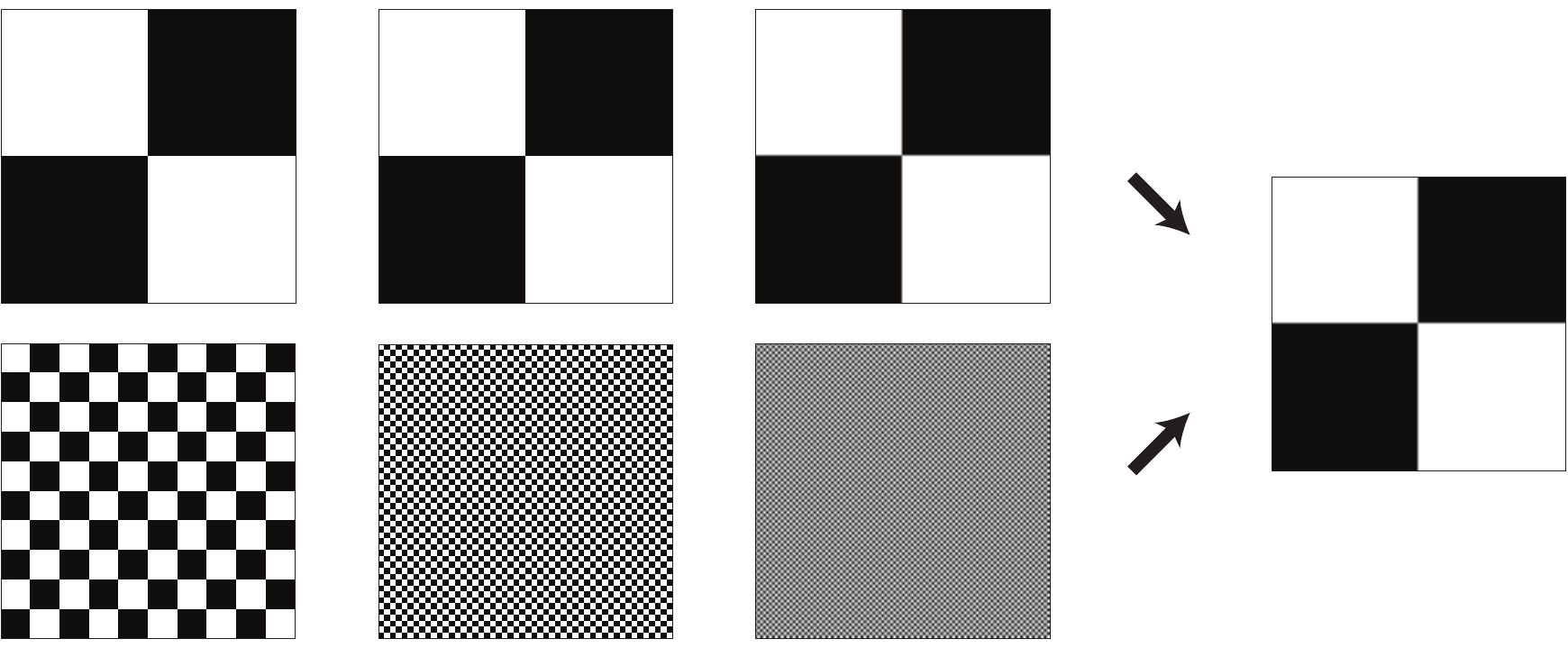}\end{center}
\vskip5pt

Homomorphism densities extend naturally to graphons. For a finite graph $G$, the density $t(\twog,G)$ may be computed by giving each vertex of $G$ a mass of $1/n$ and integrating the edge indicator function over all pairs of vertices.  In exactly the same way, the edge density $t(\twog,W)$ of a labeled graphon $W$ is
\[\int_{[0,1]^2} W(x,y) \ dx dy,\]
and the 4-cycle density $t(\fourg,W)$ is
\begin{align*}
\int_{[0,1]^4} &W(x_1,x_2)W(x_2,x_3) \quad \\
& \quad W(x_3,x_4)W(x_4,x_1) \ dx_1 dx_2 dx_3 dx_4.
\end{align*}
It is straightforward from here to write down the expression for the homomorphism density $t(H,W)$ of a finite graph $H$ into a graphon $W$. This allows us to see how the constant graphon $W \equiv 1/2$ solves the minimization problem: $t(\twog,W) = 1/2$ while $t(\fourg,W) = 1/16$.

To see the space of graphons as the completion of the space of finite graphs and make graphon convergence precise, define the \emph{cut distance} $\delta_\square (W,U)$ between two labeled graphons $W$ and $U$ by
\begin{align*}
\inf_{\varphi, \psi} \ \sup_{S, T} \bigg| \ \int \limits_{S \times T} & W\big(\varphi (x), \varphi (y)\big) \quad \\
& \quad - U \big(\psi(x),\psi(y) \big) \ dx dy \ \bigg|,
\end{align*}
where the infimum is taken over all re-labelings $\varphi$ of $W$ and $\psi$ of $U$, and the supremum is taken over all measurable subsets $S$ and $T$ of $[0,1]$. The cut distance first measures the maximum discrepancy between the integrals of two labeled graphons over measurable boxes (hence the $\square$) of $[0,1]^2$, then minimizes that maximum discrepancy over all possible re-labelings. (It is possible to define the cut distance between two finite graphs combinatorially, without any analysis, but the definition is quite involved.)

The infimum in the definition of the cut distance makes it well defined on the space of unlabeled graphons, but it is not yet a metric. Graphons $W$ and $U$ for which $t(H,W) = t(H,U)$ for all finite graphs $H$ are called \emph{weakly isomorphic}; it turns out that $W$ and $U$ are weakly isomorphic if and only if $\delta_\square (W,U) = 0$. The cut distance becomes a genuine metric on the space $\mathcal{G}$ of unlabeled graphons up to weak isomorphism. The examples of pixel picture convergence above provide examples of convergent sequences and their limits in $\mathcal{G}$ (up to weak isomorphism).

We conclude by highlighting some fundamental results on graphons.

\noindent \textbf{Theorem 1} \quad \emph{Every graphon is the $\delta_\square$-limit of a sequence of finite graphs.}

To approximate a labeled graphon $W$ by a finite labeled graph, let $S$ be a set of $n$ randomly chosen points from $[0,1]$, then construct a graph on $S$ where the edge $\{s_i,s_j\}$ is included with probability $W(s_i,s_j)$. With high probability (as $|S| \to \infty$), this labeled graph approximates $W$ well in cut distance.

\noindent \textbf{Theorem 2} \quad \emph{The space $(\mathcal{G},\delta_\square)$ is compact.}

This implies that $\mathcal{G}$ is complete; combining this fact with Theorem 1, we see that the space of graphons is the completion of the space of finite graphs with the cut metric! This theorem also demonstrates how graphons provide a bridge between different forms of Szemer\'{e}di's Regularity Lemma: Theorem 2 may be deduced from a weak form of the lemma, while a stronger regularity lemma follows from the compactness of $\mathcal{G}$. 

\noindent \textbf{Theorem 3} \quad \emph{For every finite graph $H$, the map $t(H,\cdot): \mathcal{G} \to [0,1]$ is Lipschitz continuous.}

Theorems 2 and 3 combine with elementary analysis to show that minimization problems in extremal graph theory (such as the one considered above) are guaranteed to have solutions in the space of graphons. These graphon solutions provide a ``template'', via Theorem 1, for approximate solutions in the space of finite graphs.

The interested reader is encouraged to consult L. Lov\'{a}sz's book \cite{lovasz} for more!

\bibliographystyle{plain}
\bibliography{graphonbib}

\end{multicols}

\end{document}